\documentclass[12pt, a4paper]{article}
\usepackage[left=1in,right=1in,top=1.0in,bottom=1in]{geometry}
\usepackage{amsmath,amssymb,amsthm,amsfonts,afterpage,hyperref}
\usepackage{amscd,subeqnarray}
\usepackage{graphicx}
\usepackage{caption}
\usepackage{subcaption}
\usepackage{authblk}
\usepackage{tikz}
\usepackage{array}
\usepackage{wrapfig}
\usepackage{color}
\usepackage{ragged2e}
\usepackage{stmaryrd}
\usepackage{siunitx}
\usepackage{multirow}
\usepackage{xcolor,cancel}
\usepackage{boldfonts}
\usepackage{symbols}
\usepackage{footmisc}
\usepackage{float}
\usepackage{setspace}
\usepackage{bm}
\usepackage{booktabs}
\usepackage{tikz, xcolor}
\usepackage{soul}
\usepackage[T1]{fontenc}
\usepackage[utf8]{inputenc}
\usepackage[sort, numbers]{natbib}
\setcitestyle{square}
\usetikzlibrary{shapes,arrows, positioning,calc}	
\usetikzlibrary{arrows,snakes,backgrounds}
\usepackage{algorithm}
\usepackage{algpseudocode}
\newcommand{\bfa}[1]{\boldsymbol{#1}} 			
			
\newcommand{\bfeps}{\boldsymbol{\epsilon}}

\newcommand{\Sym}{\text{Sym}}   			%
\newcommand{\curl}{\text{curl}}   				%
\newcommand{\tr}{\text{tr}}       				%

\DeclareMathAlphabet{\mathpzc}{OT1}{pzc}{m}{it}

\newcommand{\bfu}{\boldsymbol{u}}

\newcommand{\bfF}{\boldsymbol{F}}	
\newcommand{\bfE}{\boldsymbol{E}}	
\newcommand{\bfC}{\boldsymbol{C}}
\newcommand{\bfB}{\boldsymbol{B}}	
\newcommand{\bfx}{\boldsymbol{x}}	
\newcommand{\bfX}{\boldsymbol{X}}	
	
\newcommand{\bfT}{\boldsymbol{T}}		
\newcommand{\bfI}{\boldsymbol{I}}	 
\newcommand{\bfzero}{\boldsymbol{0}}

\newtheorem{theorem}{Theorem}[section]

\newtheorem{lemma}[theorem]{Lemma}

\newtheorem{Theorem}{Theorem}[section]
\newtheorem*{cf}{Continuous formulation}

\newtheorem*{dwf}{Discrete weak formulation}

\providecommand{\keywords}[1]
{
  \small	
  \textbf{\textit{Keywords---}} #1
}

\title{Finite element modeling of V-notched thermoelastic strain-limiting solids containing inclusions}

\author[1]{G. Shylaja}
\author[1]{V. Kesavulu Naidu}
\author[1]{B. Venkatesh}
\author[2]{S. M. Mallikarjunaiah}

\affil[1]{Department of Mathematics, Amrita School of Engineering, Amrita Vishwa Vidyapeetham, Bengaluru, 560035, INDIA}
\affil[2]{Department of Mathematics \& Statistics, Texas A\&M University-Corpus Christi, Corpus Christi, TX 78412, USA}

\date{} 

\begin{document}

\maketitle

\thispagestyle{empty} 

\vspace{-2em} 
\begin{center}
\small
\textsuperscript{1}g\_shylaja@blr.amrita.edu \\
\textsuperscript{1}v\_kesavulu@blr.amrita.edu (Corresponding Author) \\
\textsuperscript{1}b\_venkatesh@blr.amrita.edu  \\
\textsuperscript{2}m.muddamallappa@tamucc.edu (Corresponding Author)
\end{center}
\begin{abstract}
A precise domain triangulation is recognized as indispensable for the accurate numerical approximation of differential operators within collocation methods, leading to a substantial reduction in discretization errors. Geometrical inaccuracies are frequently introduced by conventional methods, which typically approximate curved boundaries, characteristic of thermoelastic bodies, using polygons. This often compromises the ultimate numerical accuracy. To counteract these geometrical deficiencies, methods such as isoparametric, subparametric, and isogeometric techniques have been developed, allowing for improved representation of curved surfaces. An efficient finite element method (FEM) is presented in this paper, meticulously developed to solve a complex mathematical model. This model governs the behavior of thermoelastic solids containing both a V-notch and inclusions.  The system of partial differential equations underlying this model consists of two primary components: a linear elliptic equation, which is used to describe the temperature distribution, and a quasilinear equation, which governs the mechanical behavior of the body. Through the application of this specifically tailored FEM, accurate and efficient solutions are able to be obtained for these intricate thermoelastic problems. The algebraically nonlinear constitutive equation, alongside the balance of linear momentum, is effectively reduced to a second-order quasi-linear elliptic partial differential equation. Complex curved boundaries are represented through the application of a smooth, distinctive point transformation. Furthermore, higher-order shape functions are employed to ensure the accurate computation of entries within the FEM matrices and vectors, from which a highly precise approximate solution to the BVP is subsequently obtained. The inherent nonlinearities in the governing differential equation are addressed through the implementation of a Picard-type linearization scheme. Numerical results, derived from a series of test cases, have consistently demonstrated a significant enhancement in accuracy, a crucial achievement for the nuanced analysis of thermoelastic solids.
\end{abstract}
\noindent \keywords{ Finite element method; Higher-order curved elements; Point transformation; Thermoelastic material; Strain-limiting elastic body; Quasi-linear elliptic partial differential equation.}
\section{Introduction}
The practical importance of investigating cracks and V-notches, especially within the context of a thermoelastic body, is widely acknowledged across numerous engineering disciplines. These features are fundamentally understood to act as inherent stress concentrators; consequently, in thermoelastic materials, where additional stresses are induced by temperature variations, such features become critical sites for the initiation and propagation of cracks. An understanding of their intricate behavior is therefore deemed vital for predicting component lifespan, averting catastrophic failures in essential infrastructure such as bridges, aircraft, and pressure vessels, and ultimately for safeguarding public safety. Moreover, through the accurate numerical modeling of the thermoelastic response of components containing these defects, designs can be meticulously optimized to minimize localized stress concentrations, significantly improve fatigue life, and facilitate the judicious selection of appropriate materials. This, in turn, directly contributes to the development of lighter, more efficient, and more durable products, leading to appreciable reductions in both material waste and manufacturing costs. Many contemporary technologies, including aerospace components, nuclear reactors, and high-speed machinery, are operated under extreme thermal and mechanical loads, frequently involving thermoelastic materials known to be susceptible to cracking and notching; thus, research in this domain directly contributes to the development of robust materials and designs capable of reliable performance in such demanding environments. Furthermore, insights into the thermoelastic response around cracks and notches are instrumental in the development and refinement of non-destructive testing techniques, such as thermoelastic stress analysis. These techniques enable defects to be detected and characterized without causing damage to the component, thereby facilitating timely maintenance and preventing in-service failures. Additionally, notches, particularly V-notches, are purposefully introduced in test specimens for the precise characterization of material properties, including fracture toughness and fatigue crack growth rates, under controlled stress concentration conditions; this obtained data is considered indispensable for material selection and qualification in diverse engineering projects. In essence, the rigorous study of cracks and V-notches in thermoelastic bodies is foundational to reliable engineering design, advanced material development, and the safe and efficient operation of complex systems exposed to varying thermal and mechanical conditions.

Real-world phenomena are frequently characterized by linear or nonlinear partial differential equations (PDEs), which are often defined on complex geometric domains \cite{naidu2013advantages,sasikala2023efficient}. Simulating such scenarios presents numerous challenges, primarily due to the nonlinearity of the differential operators \cite{vasilyeva2023,gou2023computational,gou2023finite}, the intricate geometric complexities of the domains \cite{MVSMM2023}, and the nonlocality inherent in the boundary conditions \cite{ferguson2015numerical}. In many instances, the differential operators are approximated using standard collocation methods, such as finite difference methods, finite element methods, boundary element methods, spectral methods, or, in some cases, meshless methods. For these techniques, the development of a high-quality computational mesh is notably challenging, yet it is essential for enhancing approximation accuracy. For example, in clinical research, images generated by computer-aided design software must be accurately represented by meshes (also known as triangulation). In the context of finite element computations, triangulation composed of straight edges (e.g., triangles or quadrilaterals) offers a straightforward approach for executing two-dimensional or three-dimensional numerical integrations through the application of highly accurate quadrature rules. This is primarily attributed to the availability of standard shape functions (such as Lagrange-type functions or other standard variants) and the accessibility of quadrature rules on regular polygons.

However, this straightforward methodology becomes inapplicable when the domain of interest possesses curved boundaries. The accurate triangulation of such a domain necessitates the incorporation of elements, typically triangles or quadrilaterals, featuring one curved edge and two straight edges to precisely represent the computational domain, thereby minimizing, or even achieving zero, discretization error. The significance of this issue is further escalated when curved inclusions or heterogeneities exist within the designated computational area. Consequently, a conceptual framework for creating curved elements was proposed \cite{ergatoudis1968curved,mcleod1975use,mitchell1979advantages,zienkiewicz2005finite}, with the accuracy of the approximation of the differential operators in such curved finite element methods (CFEM) being directly contingent upon the fidelity of the geometrical representation \cite{ciarlet1972interpolation,scott1973finite}. CFEM encompasses three principal approaches: first, the transformation of the entire global domain, including curved boundary segments, into a standard polygonal form to facilitate a global computing methodology; second, the triangulation of the curved region into a regular domain; and third, the employment of isoparametric finite elements. Each of these methods is accompanied by its own inherent challenges \cite{rathod2008use,nagaraja2010use,naidu2010use,naidu2013advantages}.

In this paper, the utilization of parabolic arcs is investigated to align with the curved boundaries (specifically, the circular inclusions within the domain) through unique point transformations, thereby leading to the derivation of a higher-order finite element method for cubic-order curved triangular elements. A numerical technique is examined for approximating the solution to a quasi-linear partial differential equation. The comprehensive algorithm developed for this purpose incorporates Picard's linearization to address the nonlinearities and utilizes a specific finite element method for spatial discretization. The mathematical model examined in this study governs the response of a limiting-strain elastic solid, as formulated within an innovative theory \cite{rajagopal2007elasticity,rajagopal2003implicit,rajagopal2011non,rajagopal2014nonlinear,rajagopal2009class,rajagopal2011modeling,gou2015modeling,MalliPhD2015,bulivcek2015analysis,bulivcek2015existence,bulivcek2014elastic}. This framework provides an appealing approach for modeling the response of a broad class of materials subjected to both mechanical and thermal loading \cite{yoon2022CNSNS, bonito2020finite}. It is well-documented that certain materials exhibit nonlinear behavior even when subjected to strains of less than 2\%, thereby underscoring the necessity for developing first-order yet nonlinear constitutive relations. Furthermore, a numerical tool capable of accurately approximating the curved regions within the computational domain is considered critical, and this paper represents the initial endeavor in this particular direction.

A particularly powerful aspect of Rajagopal's theory is its capacity to yield a hierarchy of 'explicit' nonlinear relationships, through which linearized strain is permitted to be expressed as a nonlinear function of stress. Of critical importance is the unique 'strain-limiting' property exhibited by a specific subclass of these implicit models: linearized strain is capable of being represented as a uniformly bounded function throughout the entire material domain, even when the material is subjected to conditions of significant or extreme stress. This inherent characteristic renders these models exceptionally well-suited for the investigation of crack and fracture behavior in brittle materials \cite{rajagopal2011modeling,gou2015modeling,mallikarjunaiah2015direct,MalliPhD2015}, thereby offering a robust pathway for the comprehensive analysis of both quasi-static and dynamic crack evolution, which is often problematic for conventional linear elastic fracture mechanics. The considerable utility of these strain-limiting models has been robustly demonstrated through numerous prior studies, wherein classical elasticity problems have been revisited and new, profound insights have been provided \cite{kulvait2013,rajagopal2018note,bustamante2009some,bulivcek2014elastic,itou2018states,itou2017nonlinear,yoon2022CNSNS,yoon2022MMS,rodriguez2021stretch}. A significant advantage is further afforded by their notable versatility in elucidating the mechanical behavior of a broad spectrum of materials, particularly concerning complex crack and fracture phenomena. For instance, it has been shown by recent research that when quasi-static crack evolution problems are formulated within this strain-limiting framework, the prediction of complex crack patterns and even increased crack-tip propagation velocities can be achieved \cite{lee2022finite, yoon2021quasi}. This highlights the predictive power and practical relevance of these models in advanced fracture mechanics.

This paper is meticulously organized to systematically present the theoretical underpinnings, formal problem formulation, practical numerical implementation, and comprehensive results of this investigation into the thermoelastic response of strain-limiting materials. In Section~\ref{math_formulation}, we introduce the foundational implicit theory that governs the material behavior, followed by a detailed derivation of the specific nonlinear constitutive relation utilized throughout the study. The subsequent Section~\ref{bvp_existence} is dedicated to the formal establishment of a well-posed boundary value problem for the thermoelastic material, crucially incorporating the newly derived nonlinear constitutive relation. Within this same section, we rigorously present the detailed variational (weak) formulation of the problem, accompanied by a theorem that confirms the existence of its weak solution, thereby ensuring mathematical rigor. Our numerical approach is thoroughly elucidated in Section~\ref{fem}, detailing the discretization strategy employing continuous Galerkin-type finite elements and their integration with Picard's iterative algorithm for effectively handling the inherent material nonlinearity. A comprehensive presentation and in-depth discussion of the obtained numerical solutions, along with an analysis of key parameter influences, are provided in Section~\ref{num_exp}. Finally, the paper concludes by synthesizing the primary findings and their broader implications, offering insights into their significance, and outlining compelling prospective directions for future research.
\section{Mathematical Formulation} \label{math_formulation}
This section is dedicated to providing a concise introduction to an implicit nonlinear theory, which was formulated in a series of papers by Rajagopal \cite{rajagopal2003implicit,rajagopal2007elasticity,rajagopal2011non,rajagopal2011conspectus,rajagopal2014nonlinear,rajagopal2007response} to describe the response of bulk materials. Although the constitutive relationships established within this framework are algebraically nonlinear, a mathematical model for the response of an isotropic linear elastic body, inclusive of both a V-notch and inclusions, is developed utilizing these very relations. The ultimate aim of this investigation is centered on the provision of a convergent numerical method for the accurate discretization of the resultant BVP.

Let the elastic body under consideration be denoted by $\mathcal{B}$, which is assumed to be in equilibrium under externally applied mechanical loading, with its boundary denoted by $\partial\mathcal{B}$. The elastic material is posited to occupy a region within a two-dimensional space. Points in the reference and deformed configurations of the body are denoted by $\bfX$ and $\bfx$, respectively. The linear space of symmetric tensors, equipped with the standard Frobenius norm, is denoted by $\Sym$. The displacement field is represented by $\bfu \colon \mathcal{B} \to \mathbb{R}^2$, the Cauchy stress by $\bfT \colon \mathcal{B} \to \mathbb{R}^{2 \times 2}_{\Sym}$, and the infinitesimal strain tensor by $\bfeps \colon \mathcal{B} \to \mathbb{R}^{2 \times 2}_{\Sym}$. Furthermore, the deformation gradient is denoted by $\bfF \colon \mathcal{B} \to \mathbb{R}^{2 \times 2}$, the right Cauchy-Green stretch tensor by $\bfC \colon \mathcal{B} \to \mathbb{R}^{2 \times 2}$, the left Cauchy-Green stretch tensor by $\bfB \colon \mathcal{B} \to \mathbb{R}^{2 \times 2}$, and the Lagrange strain by $\bfE \colon \mathcal{B} \to \mathbb{R}^{2 \times 2}$. 

The constitutive relationships of Cauchy elasticity were generalized by Rajagopal \cite{rajagopal2007elasticity} through the introduction of implicit elastic response relations, expressed as:
\begin{equation}\label{implicit-1}
 \mathcal{F}(\bfT, \; \bfB ) = \bfzero
\end{equation}
Here, $\mathcal{F}$ is a tensor-valued, isotropic function. A special subclass of Equation~\eqref{implicit-1} was considered by Rajagopal \cite{rajagopal2007elasticity}:
\begin{equation}\label{SL1}
\bfB := \mathcal{F}( \bfT), \quad \mbox{with} \quad \sup_{\bfT \in \Sym} \| \mathcal{F}( \bfT) \| \leq M, \;\; M >0.
\end{equation}
If the existence of such a constant $M$ is established, then the response relations given by Equation~\eqref{SL1} are referred to as \textit{strain-limiting} \cite{MalliPhD2015,mallikarjunaiah2015direct}.

Building upon the foundational work presented in \cite{itou2018states}, we now delineate several crucial properties of the tensor-valued function \eqref{SL1}, as formally stated in the following lemma. These properties are essential for understanding the behavior and well-posedness of the mathematical model.

\begin{lemma}\label{lemma1}
Let there exist positive constants $M$ and $C$ such that for all symmetric tensors $\bm{T}_1, \bm{T}_2 \in \text{Sym}(\mathbb{R}^{2 \times 2})$, the following conditions are satisfied:
\begin{itemize}
\item \textbf{Uniform Boundedness:} The norm of the tensor-valued function $\mathcal{F}(\bm{T})$ is uniformly bounded by a constant $M$. This ensures that the function's output does not grow indefinitely, regardless of the input tensor, expressed as $$\| \mathcal{F} (\bm{T}) \| \leq M.$$
\item \textbf{Monotonicity:} The function $\mathcal{F}(\bm{T})$ exhibits monotonicity. This property is crucial for ensuring the uniqueness and stability of solutions in many physical and mathematical problems. It is mathematically expressed as the inner product of the difference in function values and the difference in input tensors being non-negative: 
\[
 \left( \mathcal{F} (\bm{T}_1) - \mathcal{F} (\bm{T}_2) \right) \colon \left( \bm{T}_1 - \bm{T}_2 \right) \geq 0.
 \]
\item \textbf{Lipschitz Continuity:} The function $\mathcal{F}(\bm{T})$ satisfies Lipschitz continuity. This property provides a bound on how rapidly the function can change, which is vital for numerical stability and convergence analyses. The constant $C$ quantifies this continuity. This condition is formulated as: 
\[
\quad \dfrac{\left( \mathcal{F} (\bm{T}_1) - \mathcal{F} (\bm{T}_2) \right) \colon \left( \bm{T}_1 - \bm{T}_2 \right) }{\| \bm{T}_1 - \bm{T}_2 \|^2} \leq C.
\]
\end{itemize}
\end{lemma}
The conditions stipulated in {Lemma 2.1} collectively guarantee that $\mathcal{F}(\bm{T})$ as a {monotone operator} across its entire domain, $\text{Sym}(\mathbb{R}^{2 \times 2})$. This characteristic is of paramount importance in the mathematical analysis of nonlinear systems, as it often underpins the existence, uniqueness, and stability of solutions to associated boundary value problems.

Upon the application of the standard linearization procedure, the following form for $\bfeps$ is obtained:
\begin{equation}
\bfeps = \beta_1 \, \bfI + \beta_2 \, \bfT + \beta_2 \, \bfT^2.
\end{equation}
Finally, the governing system of equations, utilized to model the behavior of an elastic material within the framework of algebraically nonlinear theories, is given by:
\begin{subequations}
\begin{align}
-\nabla \cdot \bfT &=\bf{0}, \quad \mbox{and} \quad \bfT = \bfT^\mathrm{T}, \label{equilib:eq} \\
\bfeps &= \Psi_{0}\left( \tr \, \bfT, \; \| \bfT \| \right) \bfI + \Psi_{1}\left( \| \bfT \| \right) \bfT, \quad {\Psi_{0}\left( 0, \; \cdot \right)=0}, \label{eqn:main} \\
\curl \, \curl \, \bfeps &=\bfa{0}, \label{eqn_scompata} \\
\bfeps &:= \frac{1}{2} \left( \nabla \bfu + \nabla \bfu^{T} \right). \label{eqn_linstrain}
\end{align}
\end{subequations}
Both $\Psi_{0} \colon \mathbb{R} \times \mathbb{R}_{+} \to \mathbb{R}$ and $\Psi_{1} \colon \mathbb{R}_{+} \to \mathbb{R}$ are functions of the invariants of the Cauchy stress. In the above governing system of equations, the first component of Equation~\eqref{equilib:eq} is derived from the balance of linear momentum for a quasi-static situation, where negligible body forces are assumed to be acting on the material. The second component of Equation~\eqref{equilib:eq} is obtained from the balance of angular momentum. Equation~\eqref{eqn:main} represents the constitutive relationship that governs the material's response to both mechanical and thermal stimuli. Furthermore, the strain-compatibility condition is expressed by Equation~\eqref{eqn_scompata}, and the equation for the linearized strain is given by Equation~\eqref{eqn_linstrain}.

In the context of anti-plane shear deformation, the displacement field, denoted as $\bm{u}$, simplifies significantly. In this particular configuration, the displacement is exclusively a scalar function of the spatial coordinates $x$ and $y$, such that the displacement vector is expressed as:
\begin{equation}\label{eq:disp_vector}
\bm{u}(x,\, y) = \left( 0, 0, w( x,\, y) \right).
\end{equation}
This implies that deformation occurs only out of the $xy$-plane, along the $z$-axis.

For a linear, isotropic, and homogeneous material, the classical constitutive relationship between the stress tensor $\bm{T}$ and the infinitesimal strain tensor $\bm{\epsilon}$ is given by:
\begin{equation}
\bm{T} = 2 \, \mu \, \bm{\epsilon},
\end{equation}
where $\mu$ represents the shear modulus of the material, a crucial material property that quantifies its resistance to shear deformation, with units of N/m$^2$. This equation establishes a direct linear proportionality between stress and strain.

However, within the framework of the \textit{strain-limiting theory of elasticity}, a more generalized constitutive equation is employed. This theory accounts for situations where the material's response may deviate from linearity, especially at higher strain levels. In this advanced theory, the strain is expressed as a function of the stress magnitude:
\begin{equation}\label{eqn:main2_strain}
\bm{\epsilon} = \Psi\left( \| \bm{T} \| \right) \bm{T}.
\end{equation}
This relationship implies that the material's compliance is dependent on the magnitude of the applied stress. It is crucial to note that this equation is assumed to be invertible \cite{mai2015strong,mai2015monotonicity,rajagopal2007response,rajagopal2011conspectus,rajagopal2004thermomechanical}, a property supported by various studies in the field. The invertibility of this relation signifies that a unique stress state corresponds to a given strain state.

Consequently, the inverted constitutive equation, which is characteristic of hyperelastic materials, expresses the stress components in terms of the strain magnitude:
\begin{equation}\label{eqn:main2_stress}
\bm{T} = \Phi\left( \| \bm{\epsilon} \| \right) \bm{\epsilon}.
\end{equation}
This formulation highlights the hyperelastic nature of the material within this theory, where the stress can be derived from a strain energy potential.

Finally, the principle of balance of linear momentum, when applied to this constitutive relationship, yields a quasilinear second-order PDE. This governing equation, which describes the equilibrium of forces within the deformed body, takes the form:
\begin{equation}\label{pde:nlin1}
- \nabla \cdot \left( \Phi \left( \| \bm{\epsilon} \| \right) \;  \bm{\epsilon} \right) =0.
\end{equation}
This non-linear PDE is a cornerstone for analyzing the deformation and stress distribution in materials exhibiting strain-limiting behavior under anti-plane shear conditions. For the purpose of developing the BVP to study the state of response of thermoelastic body, we consider the following form for 
$$\Phi \left( s \right) = \dfrac{1}{ \left( 1 - (\beta \, s)^\alpha\right)^{1/\alpha}}$$


\section{BVP for Thermoelastic State and Existence of Solution}\label{bvp_existence}
This section focuses on the mathematical formulation and conditions for the existence of solutions for the problem of the response of a thermoelastic solid containing voids and inclusions. We're dealing with a solid body whose behavior is described by an {algebraically nonlinear relationship}, even though the material itself is considered {geometrically linear}. This means that while the relationship between stress and strain isn't strictly proportional, the deformations are assumed to be small enough that the geometry of the body doesn't significantly change. The material is also assumed to be {homogeneous} and {initially unstressed}.

The solid body occupies a domain $\mathcal{B} = \Omega \times \mathbb{R}$, where $\Omega$ is a simply connected domain in $\mathbb{R}^2$. Its boundary, $\partial \Omega$, is assumed to be sufficiently smooth ($C^{0,\,1}$). This boundary is divided into two distinct, non-overlapping parts: $\Gamma_D$ and $\Gamma_N$. These represent areas where {Dirichlet (displacement) boundary conditions} and {Neumann (traction) boundary conditions} are applied, respectively. Each of these boundary segments can consist of a finite number of individual parts, and each can be parameterized for mathematical convenience.

To analyze this problem rigorously, we introduce several function spaces. The space of \textit{Lebesgue integrable functions} is denoted by $L^{p}(\Omega)$, where $p\in[1, \infty)$. For $p=2$, we use $L^2(\Omega)$ to represent the space of square-integrable functions. The inner product in $L^2(\Omega)$ is denoted by $\left( \cdot, \; \cdot \right)_{L^2}$, and the corresponding norm is $\| \cdot \|_{L^2}$. The \textit{classical Sobolev space}, $H^{1}(\Omega)$, plays a crucial role in the analysis of partial differential equations. Its norm is defined as:
\[
\| v \|^2_{H^{1}(\Omega)} := \| v \|^2_{L^2} + \| \nabla v \|^2_{L^2},
\]
where the terms on the right-hand side refer to the $L^2$ norm. Additionally, $H^{-1}(\Omega)$ is defined as the dual space to $H_0^{1}(\Omega)$, which is a subspace of $H^1(\Omega)$ containing functions that vanish on the entire boundary. We also define two important subspaces of $H^1(\Omega)$:
\begin{equation}
V := H^1(\Omega), \quad V^{0} := \left\{ v \in H^{1}(\Omega) \colon v=0 \; \mbox{on} \; \Gamma_D\right\}.
\end{equation}

In this particular problem, the out-of-plane displacement function, $w(x, y)$, is the sole unknown variable that fully describes the mechanical behavior of the material body. For understanding the state of temperature distribution in the body we couple the mechanics equation with the temperature variable $\theta(x, \, y)$.  Due to the specific loading conditions, only two components in the stress and strain tensors are non-zero, significantly simplifying the problem.

The problem, as considered here, is formally expressed as:
\begin{cf}
Given all the parameters in the model, we seek to find the functions $\theta \in C^2(\Omega)$ and $w \in C^2(\Omega)$ such that:
\begin{subequations}\label{f1a}
\begin{align}
- \nabla \cdot \left( \kappa \, \nabla \theta \right) &= g, \quad \mbox{on} \quad  \Omega \label{f1_1} \\
- \nabla \cdot \left( \Phi \left( \|  \bm{\epsilon} \| \right) \;  \bm{\epsilon} \right)  + \xi \, \theta &=0. \quad \mbox{on} \quad  \Omega \label{f1_2}
\end{align}
\end{subequations}
\end{cf}
In this reformulation, the original problem, which might have involved both traction and displacement conditions, has been transformed into a {non-homogeneous Dirichlet-type problem} for a single unknown function, $w$. This approach is consistent with methodologies explored in previous research, as seen in works like \cite{bulivcek2014elastic, bulivcek2015existence, kulvait2019state, kulvait2013, yoon2022MMS}.

\subsection{Variational Formulation and Existence of Solution}
This subsection develops a {variational formulation} for the nonlinear boundary value problem previously discussed. To achieve this, we first lay out several crucial assumptions regarding the problem's input data. For the robust formulation and solvability of the nonlinear boundary value problem, the following critical assumptions are imposed on the material properties and relevant data:
\begin{itemize}
\item[A1] \textbf{Thermal Conductivity ($k$)}: The heat conduction coefficient, $k$, is strictly positive and bounded. For a homogeneous material, $k$ is a constant within $\Omega$ such that $0 < k_1 < k < k_2 < \infty$. For a heterogeneous material, $k(\bm{x})$ is an $L^\infty$ function from $\Omega$ to $\text{Sym}(\mathbb{R}^{d \times d})$, with its essential infimum ($k_1$) and essential supremum ($k_2$) satisfying:
\[
0 < k_1 := \operatorname{essinf}_{\bm{x} \in \Omega} \inf_{\bm{v} \in \mathbb{R}^d, \, \bm{v} \neq \bm{0}} \dfrac{(k(\bm{x})\bm{v}, \, \bm{v} )}{(\bm{v}, \, \bm{v})}, \quad \infty > k_2 := \operatorname{esssup}_{\bm{x} \in \Omega} \sup_{\bm{v} \in \mathbb{R}^d, \, \bm{v} \neq \bm{0}} \dfrac{(k(\bm{x})\bm{v}, \, \bm{v} )}{(\bm{v}, \, \bm{v})}.
\]
\item[A2] \textbf{Elastic and Thermal Material Coupling Constant ($\mu, \xi$)}: These constants are assumed to be bounded and positive. For homogeneous materials, $\mu, \xi$ are constant. For heterogeneous materials, they are $L^\infty$ functions on $\Omega$, i.e., $\mu \in L^{\infty}(\Omega, \, \mathbb{R})$. Specifically for the shear modulus $\mu(\bm{x})$, we require:
\[
0 < \mu_1 := \operatorname{essinf}_{\bm{x} \in \Omega} \mu(\bm{x}) \leq \operatorname{esssup}_{\bm{x} \in \Omega} \mu(\bm{x}) =: \mu_2 < \infty.
\]
\item[A3] \textbf{Data Regularity}: The body force $g$, initial displacement $w_0$, and initial temperature $\theta_0$ are required to possess sufficient regularity. Specifically, $g \in L^{2}(\Omega)$, $w_0 \in \left( L^2(\Omega) \right)^2$, and $\theta_0 \in L^2(\Omega)$.
\end{itemize}
Extensive mathematical proofs for the existence and uniqueness of weak solutions to problems involving strain-limiting constitutive descriptions have been established in previous works, notably \cite{bulivcek2014elastic,bulivcek2015existence,beck2017existence,bulivcek2015analysis,bonito2020finite}. In this paper, we aim to reframe the current problem to directly enable a comparative analysis with these prior studies, particularly with \cite{beck2017existence}, concerning the existence and uniqueness of its solution within the context of strain-limiting theory.

\subsection*{Existence of Solution Theorem}

Consider a bounded, connected, and Lipschitz domain $\Omega \subset \mathbb{R}^d$. Its boundary, $\partial \Omega$, is composed of two disjoint open sets: a Dirichlet boundary $\Gamma_D$ and a Neumann boundary $\Gamma_N$, such that their closures cover the entire boundary ($\overline{\Gamma_D \cup \Gamma_N} = \partial \Omega$).

\begin{Theorem}
Let $\Omega \subset \mathbb{R}^2$ be bounded, connected, Lipschitz domain with open sets Dirichlet boundary $\Gamma_D$ and Neumann boundary $\Gamma_N$ such that $\Gamma_D \, \cap \, \Gamma_N = \emptyset$ and $\overline{\Gamma_D \cup \Gamma_N} = \partial \Omega$. Given a vector field $\bm{f} \colon \Omega \to \mathbb{R}^2$, a traction $\bm{g}: \Gamma_N \to \mathbb{R}^2$, a boundary displacement $\bm{u}_0 \colon \Gamma_D \to \mathbb{R}^2$, and a bounded mapping $\mathcal{F} \colon \text{Sym}(\mathbb{R}^{2 \times 2}) \to \text{Sym}(\mathbb{R}^{2 \times 2})$, then find the pair $(\bm{u}, \, \bm{T})$ such that $\bm{u} \colon \overline{\Omega} \to \mathbb{R}^{2}$, and $\bm{T} \colon \overline{\Omega} \to \text{Sym}(\mathbb{R}^{2 \times 2})$, and
\begin{align}
- \div \; \bm{T} &= \bm{f} \quad \mbox{in} \quad \Omega, \notag \\
\bm{\epsilon}(\bm{u}) &= \mathcal{F}(\bm{T}) := \dfrac{\bm{T}}{(1 + \left( \beta \, \| \bm{T} \| \right)^a)^{1/a}},\quad \; a>0, \,\beta \geq 0 \quad \mbox{in} \quad \Omega, \notag \\
\bm{u} &= \bm{u}_0, \quad \mbox{on} \quad \Gamma_D, \\
\bm{T} \bm{n} &= \bm{g}, \quad \mbox{on} \quad \Gamma_N \notag
\end{align}
with following assumptions on the data:\\
\noindent [B1] $\bm{f} \in L^{2}(\Omega)^2$, $\bm{g} \in L^{2}(\Gamma_N)^2$, and
\[
\bm{0} = \int_{\Omega} \bm{f} \, d\bm{x} + \int_{\partial \Omega} \bm{g} \, dS
\]
[B2] $\bm{u}_0 \in W^{1, \, \infty}(\Omega)^2$, with $\nabla \bm{u}_0(\bm{x})$ for almost every $\bm{x} \in \overline{\Omega}$ contained in a compact set in $\mathbb{R}^{2 \times 2}$.
\newline Assume that the data $(\bm{f}, \; \bm{g}, \; \bm{u}_0)$ satisfy the above assumptions [B1]-[B2], and consider $a >0$,
then there exists a pair $(\bm{u}, \, \bm{T}) \in W^{1, \; \infty}(\Omega)^{2} \times \text{Sym}(L^1(\Omega))^{2 \times 2}$ satisfying
\begin{equation}
\int_{\Omega} {\bm{T} \cdot \bm{\epsilon}(\bm{w})} \, d\bm{x} = \int_{\Omega} \bm{f} \cdot \bm{w} \, d\bm{x}, \quad \forall \; \bm{w} \in C^{1}_{0}(\Omega)^2
\end{equation}
\end{Theorem}
This theorem was originally proven in \cite{beck2017existence} for a specific representation of the "stress tensor". Our current problem formulation is equivalent to the one studied in \cite{beck2017existence} under specific conditions. Specifically, we consider the body force ${f} = -\xi  \,  \theta$ and a vanishing traction boundary condition $\bm{g} = \bm{0}$. Additionally, the Cauchy stress tensor in our formulation is given by:
\begin{equation}
\bm{T} :=\mathcal{L}(\bm{\epsilon}(\bm{u})),
\end{equation}
where $\mathcal{L}(\cdot)$ is a uniformly monotone operator with at most linear growth at infinity. Given that the material parameters $k, \, \mu, \, \xi, \, \theta_0, \, \bm{u}_0$ all satisfy the assumptions (A1)-(A3), the existence of a unique pair $(\bm{u}, \; \bm{T}) \in (H_0^1(\Omega))^2 \times \text{Sym}(L^1(\Omega)^{2 \times 2})$ is guaranteed.

Given the complexity of the robust formulation presented earlier, it's important to note that an exact analytical solution isn't available, even for simplified one-dimensional scenarios. Because of this, our focus shifts to developing a stable and convergent numerical method to approximate the solution to the formulation outlined before. In the upcoming section, we'll introduce a novel {finite element discretization} approach. What makes this method unique is its use of a point transformation during the calculation of shape functions. This specific transformation, combined with an optimized integration procedure, consistently yields superior numerical results.

\section{Finite Element Discretization of the Nonlinear Thermoelastic Model}\label{fem}
This paper introduces a novel finite element discretization approach for a nonlinear BVP. This BVP accurately models the behavior of strain-limiting elastic materials that feature both a V-notch and multiple inclusions. Accurately solving this problem numerically presents two primary challenges: achieving an accurate mesh discretization of the domain and effectively approximating the PDE. To overcome the discretization challenge, we employ curved triangular elements for triangulating the domain. A key advantage of this method is the elimination of discretization errors, as these elements are constructed using the exact geometry of the problem. For the numerical approximation of the PDE, we utilize cubic-order shape functions. These functions are derived from a unique point transformation method well-established in existing literature \cite{mcleod1975use,rathod2008use}. The effectiveness of these higher-order curved triangular elements has been extensively demonstrated in previous research \cite{shylaja2021improved, shylaja2021two, Murali2019DBF1, shylaja2019finite, Murali2019DB, Murali2019DBF2}.

Since an exact analytical solution for the specific physical setup considered in this study is unavailable, we assess the quality of our proposed approximation by calculating the relative difference in the numerical solution between successive Picard iterations. This metric provides a robust measure of convergence and accuracy.

Let $\Omega$ represent a two-dimensional open domain, with its boundary denoted by $\Gamma$. We assume that $\Gamma$ is sufficiently smooth. The boundary $\Gamma$ is composed of two distinct sets: $\Gamma_{w}$ and $\Gamma_h^i$. Specifically, $\Gamma = \Gamma_{w} \cup_{i=1}^{m} \Gamma_h^i$, where $m \in \mathbb{N}$ signifies the total number of inclusions. Here, $\Gamma_{w}$ corresponds to the portion of the boundary where Dirichlet boundary conditions are applied. Conversely, $\Gamma_h^i$ represents the internal inclusions, which are assumed to be traction-free. The regions immediately surrounding these inclusions are triangulated using specialized "curved" triangles, characterized by one curved edge and two straight edges.

To approximate the primary variables, $w$ and $\theta$, we construct a discrete solution space, $S_h$, defined as follows:
\begin{equation}
S_h = \left\{ v_h \in \left( C(\overline{\Omega})\right) \colon \left. v_h\right|_K \in \mathbb{Q}_3, \; \forall K \in \mathcal{T}_h \right\}
\end{equation}
In this formulation, $\mathbb{Q}_3$ represents a collection of cubic-order shape functions that were previously developed for the reference element, $\widehat{K}$. The final discrete approximation space, denoted as $\widehat{V}_h$, is then derived by intersecting $S_h$ with the continuous solution space, $V$. This intersection ensures that the approximated solutions not only adhere to the finite element discretization but also satisfy any global constraints or essential boundary conditions inherent to the problem. The relationship between these spaces is given by:
\begin{equation}\label{app-spaces}
\widehat{V}_h = S_h \, \cap \, V, \quad \mbox{and} \quad \widehat{V}_{h, 0} =  \left\{ \varphi_h \in \widehat{V}_h : \varphi_i = 0 \text{ for all } i \in \Gamma_D \right\}, 
\end{equation}

This meticulous definition of the approximation space is crucial for ensuring the accuracy and stability of our numerical scheme.

\subsection{Discrete Formulation and Picard's Iteration for Nonlinear Problems}
To transition from the continuous problem to a solvable form, we derive a discrete weak formulation by employing the shape functions established in the approximation space defined earlier in equation \eqref{app-spaces}. The primary challenge in this process lies in effectively managing the nonlinearities present in the continuous governing equations. These nonlinearities can be addressed either at the level of the differential equation itself or during the assembly of the linear algebraic system. In this work, we opt for the former approach, utilizing the widely recognized Picard's type linearization method. This technique systematically transforms the nonlinear problem into a sequence of tractable linear problems.

The efficiency and convergence rate of this iterative algorithm are highly dependent on a judicious choice of the initial guess. To ensure rapid convergence in our implementation, we initially solve a simplified linear problem derived from linearizing the quasilinear weak formulation. The solution obtained from this preliminary step then serves as the initial guess for the subsequent Picard iterations. Our numerical simulations consistently demonstrated convergence within a reasonable number of iterations, indicating the robustness of our approach. Overall, our discrete formulation guarantees a unique solution, which is a crucial property for reliability.

The final discrete finite element problem is formally stated as follows:
\begin{dwf}
Given all the parameters in the model and the Dirichlet boundary data for both the variables, find $\theta \in \widehat{V}_h$ such that 
\begin{equation}
\sum_{K \in T_h} \int_K  \kappa \,  \nabla \theta_h \cdot \nabla \varphi_h \,d\bfa{x}  = 0, \quad \forall \varphi_h \in \widehat{V}_{h, \, 0},   
\end{equation}
and give  $n^{th}$ iteration solution $w^n_h \in \widehat{V}_h$, for $n=0, 1, 2, \cdots $, find $w^{n+1}_h \in \widehat{V}_h$ such that
\begin{equation}
\sum_{K \in T_h} \int_K \Phi(w^n_h) \;  \nabla w^{n+1}_h \cdot \nabla \psi_h  \; d\bfa{x} = - \sum_{K \in T_h} \int_K  \xi \, \theta \, \psi \; d\bfa{x}   \quad \forall \psi_h \in \widehat{V}_{h, \, 0}. 
\end{equation}
\end{dwf}
Next, we briefly explain the derivation of the shape functions utilizing the cubic-order point transformations for the triangular elements. 
\subsection{Cubic-order Point Transformations for Triangular Elements with One Curved Boundary}
\begin{figure}[H]
\centering
\noindent
\begin{subfigure}[b]{0.4\textwidth}
\includegraphics[width=1\textwidth]{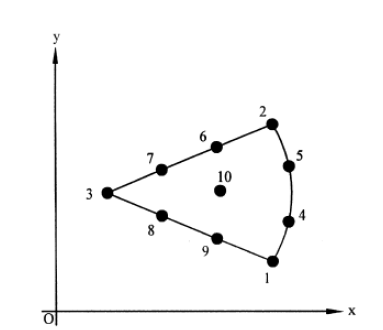}
\caption{Unmapped triangular element of cubic order}
\end{subfigure}
\qquad
\noindent
\begin{subfigure}[b]{0.4\textwidth}
\includegraphics[width=1\textwidth]{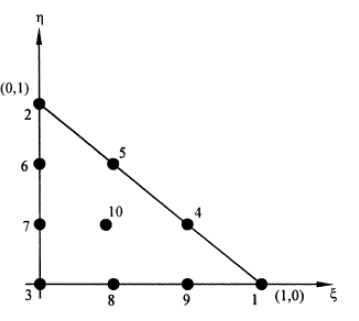}
\caption{Mapped triangular element of cubic order}
\end{subfigure}
\caption{Mapping the 10-node cubic order with one edge curved triangle to the standard right-angled unit triangle}
\label{fig_mapping}
\end{figure}

We consider specialized triangular elements where one edge exhibits a curvature while the remaining two edges are straight, as visually depicted in Figure \ref{fig_mapping}. For the approximation of a field variable, denoted as $\phi$, which governs the underlying physical problem, we employ Lagrange interpolants. These interpolants are defined by:
\begin{equation}\label{eq_phi}
\phi=\sum\limits_{i=1}^\frac{(n+1)(n+2)}{2}{ N_{i}^{(n)}(\xi ,\eta )\phi_{i}^{e}} \quad n=3 \text{ (here)}
\end{equation}
In this expression, $N_{i}^{(n)}(\xi,\eta)$ represents the standard shape functions for a triangular element of order $n$, as comprehensively derived in \cite{rathod2008use} and detailed in equation \eqref{SF}. These shape functions are evaluated at node $i$ within the local coordinate system $(\xi, \eta)$.

To facilitate computations between the physical (cartesian) and local (natural) coordinate systems, we establish transformation formulae. These transformations map the nodal values of the physical coordinates (e.g., $x$ or $y$) from the physical domain to the local coordinates within the reference element:
\begin{equation}
\label{eq_subpara}
t=\sum\limits_{i=1}^\frac{(n+1)(n+2)}{2}{ N_{i}^{(n)}(\xi ,\eta )t_{i}}
\end{equation}
For the specific geometry of our elements, the nodes positioned along the straight edges, namely edge $3-1$ and edge $3-2$ in Figure \ref{fig_mapping}, are uniformly spaced. When applying standard analytical geometry principles to subdivide these straight line segments, the general transformation formula in equation \eqref{eq_subpara} simplifies considerably. Specifically, if $t$ represents a nodal coordinate (either $x$ or $y$) of the triangular element, the transformation can be expressed as:
\begin{equation}
t(\xi,\eta)=t_3+(t_1-t_3)\xi+(t_2-t_3)\eta +a_{11}^{(n)}(t)\xi\eta+H(n-3)\sum_{i+j=n \atop i \neq j} a_{ij}^{(n)}\xi^i\eta^j, \quad (1\leq i, j \leq n-1),
\end{equation}
Here, $t_1, t_2, t_3$ denote the nodal values of the coordinate $t$ at vertices 1, 2, and 3, respectively. The term $H(n-3)$ is the well-known Heaviside step function (or unit step function), which, in this context, clarifies the applicability of the higher-order terms:
\begin{equation}
H(n-3) = \begin{cases}
0 \quad & \text{if } n<3\\
1 \quad & \text{if } n \geq 3
\end{cases}
\end{equation}
This function ensures that the summation term, which accounts for higher-order blending functions along the curved edge, is only active for elements of cubic order ($n=3$) or higher, allowing for accurate representation of the curved geometry.

\subsubsection{Derivation of Shape function}
The primary objective of shape functions is to represent the behavior of the unknown quantity being analyzed (e.g., displacement, temperature, stress) within an element. They provide a convenient mathematical framework to interpolate the values of the unknown quantity at any point within the element based on known values at specific locations, typically referred to as nodes. In practice, shape functions are defined and manipulated algebraically to derive element stiffness matrices, load vectors, and other quantities needed to solve the finite element equations. They are typically expressed in terms of local coordinates within an element and then transformed to global coordinates as needed for assembly into the overall system. Understanding and appropriately choosing shape functions are crucial for obtaining accurate and efficient solutions in finite element simulations. They form the backbone of FEM analysis, enabling the approximation of complex physical phenomena with simple mathematical constructs within individual finite elements.

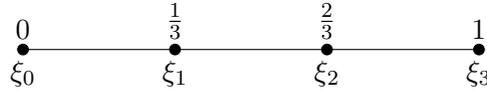
\begin{figure}[H]
	\centering
\begin{tikzpicture} [scale=2.0]
		\draw (0,0) -- (3,0);
		\filldraw [black] (0,0) circle (1pt) node[below]{\small{$\xi_0$}} node[above]{\small{0}};
		\filldraw [black] (1,0) circle (1pt) node[below]{\small{$\xi_1$}} node[above]{\small{$\frac{1}{3}$}};
		\filldraw [black] (2,0) circle (1pt) node[below]{\small{$\xi_2$}} node[above]{\small{$\frac{2}{3}$}};
		\filldraw [black] (3,0) circle (1pt) node[below]{\small{$\xi_3$}} node[above]{\small{1}};
\end{tikzpicture}
\caption{Nodal representation for derivation of shape function}
\label{fig_sf}
\end{figure}

Lagrange shape functions are fundamental in the finite element method, especially for interpolating the behavior of unknown quantities within finite elements. These shape functions are based on Lagrange interpolation, a method commonly used to approximate functions by polynomial interpolation. Lagrange interpolation is a technique for constructing a polynomial that passes through a given set of points. The Lagrange polynomial $L_i$ of degree $n$ is defined as,
\begin{equation}
L_i(x)\equiv N_i = \prod_{j=1, j\neq i}^{n} \frac{x-x_j}{x_i-x_j} \equiv \frac{x-x_1}{x_i-x_1}\times \frac{x-x_2}{x_i-x_2} \times \cdot\cdot\cdot \times \frac{x-x_n}{x_i-x_n}
\end{equation}
where $x_i$ are the interpolation nodes (also known as nodal points) and $i=0,1,2,\cdot\cdot\cdot,n$. Each Lagrange polynomial $L_i$ is associated with a specific node and has the property $L_i(x_j)=\delta_{ij}$, where $\delta_{ij}$ is the Kronecker delta.\\
For a two-dimensional finite element with $NP=\frac{(n+1)(n+2)}{2}$ nodes in the triangular element, the Lagrange shape function $N_i$ for node $i$ are defined as follows,
\begin{equation}
N_i(x) = \prod_{j=1, j\neq i}^{n} \frac{x-x_j}{x_i-x_j} 
\end{equation}
From the Figure \ref{fig_mapping}b consider the line segment from the triangular element as in Figure \ref{fig_sf}\\
Consider the following,
\begin{align}
\nonumber
    \lambda_{0}^{0} = & 1\\
    \nonumber
    \lambda_{1}^{1} = & \frac{\xi - \xi_{0}}{\xi_{1}-\xi_{0}} = \frac{\xi - 0}{1-0} = 3\xi\\ 
    \nonumber
	\lambda_{2}^{2} = & \frac{(\xi - \xi_{0})(\xi - \xi_{1})}{(\xi_{2}-\xi_{0})(\xi_{2}-\xi_{1})} = \frac{9}{2} \left(\xi^2 -\frac{\xi}{3}\right)\\
 \label{eq_lambda}
	\lambda_{2}^{2} = & \frac{(\xi - \xi_{0})(\xi - \xi_{1})(\xi - \xi_{2})}{(\xi_{3}-\xi_{0})(\xi_{3}-\xi_{1})(\xi_{3}-\xi_{2})} = \frac{1}{2} \left(9\xi^3 - 9\xi^2+2\xi\right)  
\end{align}
Each node in Figure \ref{fig_mapping} is written in terms of \eqref{eq_lambda} as below,
\begin{align}
    \nonumber
    N_{1}^{(3)}= & \lambda_{3}^{3}(\xi)\lambda_{0}^{0}(\eta)\lambda_{0}^{0}(\delta) = \frac{1}{2} \left(9\xi^3 - 9\xi^2+2\xi\right) \\
    \nonumber
    N_{2}^{(3)}= & \lambda_{0}^{0}(\xi)\lambda_{3}^{3}(\eta)\lambda_{0}^{0}(\delta) = \frac{1}{2} \left(9\eta^3 - 9\eta^2+2\eta\right)\\
    \nonumber
    N_{3}^{(3)}= & \lambda_{0}^{0}(\xi)\lambda_{0}^{0}(\eta)\lambda_{3}^{3}(\delta) = \frac{1}{2} \left(9\delta^3 - 9\delta^2+2\delta\right)\\
    \nonumber
    N_{4}^{(3)}= & \lambda_{2}^{2}(\xi)\lambda_{1}^{1}(\eta)\lambda_{0}^{0}(\delta) = \frac{9}{2} \left(\xi^2 -\frac{\xi}{3}\right)3\eta\\
    \nonumber
    N_{5}^{(3)}= & \lambda_{1}^{1}(\xi)\lambda_{2}^{2}(\eta)\lambda_{0}^{0}(\delta) = 3\xi\left[\frac{9}{2} \left(\eta^2 -\frac{\eta}{3}\right)\right]\\
    \label{eq_SFlambda}
    N_{6}^{(3)}= & \lambda_{0}^{0}(\xi)\lambda_{2}^{2}(\eta)\lambda_{1}^{1}(\delta) = \frac{9}{2} \left(\eta^2 -\frac{\eta}{3}\right)3\delta\\
    \nonumber
    N_{7}^{(3)}= & \lambda_{0}^{0}(\xi)\lambda_{1}^{1}(\eta)\lambda_{2}^{2}(\delta) = 3\eta \left[\frac{9}{2} \left(\delta^2 -\frac{\delta}{3}\right)\right]\\
    \nonumber
    N_{8}^{(3)}= & \lambda_{1}^{1}(\xi)\lambda_{0}^{0}(\eta)\lambda_{2}^{2}(\delta) = 3\xi \left[\frac{9}{2} \left(\delta^2 -\frac{\delta}{3}\right)\right]\\
    \nonumber
    N_{9}^{(3)}= & \lambda_{2}^{2}(\xi)\lambda_{0}^{0}(\eta)\lambda_{1}^{1}(\delta) = \left[\frac{9}{2} \left(\xi^2 -\frac{\xi}{3}\right)\right] 3\delta\\
    \nonumber
    N_{10}^{(3)}= & \lambda_{1}^{1}(\xi)\lambda_{1}^{1}(\eta)\lambda_{1}^{1}(\delta) = (3\xi)(3\eta)(3\delta) 
\end{align}
Replacing $\delta = 1-\xi-\eta$ in \eqref{eq_SFlambda} we get the required shape function for cubic order as follows,
\begin{subequations}\label{SF}
\begin{align}
N_{1}^{(3)}(\xi ,\,\,\eta ) &=\frac{9{{\xi }^{3}}}{2}-\frac{9{{\xi }^{2}}}{2}+\xi,  \\
N_{2}^{(3)}(\xi ,\,\,\eta ) &=\frac{9{{\eta }^{3}}}{2}-\frac{9{{\eta }^{2}}}{2}+\eta,  \\
N_{3}^{(3)}(\xi ,\,\,\eta ) &=-\frac{9{{\xi }^{3}}}{2}-\frac{9{{\eta }^{3}}}{2}-\frac{27{{\xi }^{2}}\eta }{2}-\frac{27\xi {{\eta }^{2}}}{2}+9{{\xi }^{2}}+9{{\eta }^{2}}+18\xi \eta -\frac{11\xi }{2}-\frac{11\eta }{2}+1,  \\
N_{4}^{(3)}\,(\xi ,\,\eta ) &=\frac{27{{\xi }^{2}}\eta }{2}-\frac{9\xi \eta }{2} \\
N_{5}^{(3)}\,(\xi ,\,\eta ) &=\frac{27\xi {{\eta }^{2}}}{2}-\frac{9\xi \eta }{2} \\
N_{6}^{(3)}\,(\xi ,\,\eta ) &=-\frac{27{{\eta }^{3}}}{2}-\frac{27\xi {{\eta }^{2}}}{2}+18{{\eta }^{2}}+\frac{9\xi \eta }{2}-\frac{9\eta }{2} \\
N_{7}^{(3)}\,(\xi ,\,\eta ) &=\frac{27{{\eta }^{3}}}{2}+27\xi {{\eta }^{2}}-\frac{45{{\eta }^{2}}}{2}+\frac{27{{\xi }^{2}}\eta }{2}-\frac{45\xi \eta }{2}+9\eta, \\
N_{8}^{(3)}(\xi ,\,\,\eta ) &=\frac{27{{\xi }^{3}}}{2}+27{{\xi }^{2}}\eta +\frac{27\xi {{\eta }^{2}}}{2}-\frac{45{{\xi }^{2}}}{2}-\frac{45\xi \eta }{2}+9\xi, \\
N_{9}^{(3)}\,(\xi ,\,\eta ) &=-\frac{27{{\xi }^{2}}}{2}-\frac{27{{\xi }^{2}}\eta }{2}+18{{\xi }^{2}}+\frac{9\xi \eta }{2}-\frac{9\xi }{2}, \\
N_{10}^{(3)\,}(\xi ,\,\eta ) &=-27\xi {{\eta }^{2}}-27{{\xi }^{2}}\eta +27\xi \eta, 
\end{align}
\end{subequations}
The Figure \ref{fig_mapping} depicts the domain with two straight edges and one curved edge. The subparametric transformations given below are used to convert the global coordinates of a typical triangle element in cubic order to the local coordinates. The cubic order's point transformation is as follows:
\begin{equation}
    t(\xi,\eta )= {{t}_{3}}+({{t}_{1}}-{{t}_{3}})\xi +({{t}_{2}}-{{t}_{3}})\,\eta+\frac{9}{4}[({{t}_{4}}+{{t}_{5}})-({{t}_{1}}+{{t}_{2}})]\xi \eta ,\quad  t=(x,y),
\end{equation}
with the boundary node $t_{5}$ as,
\begin{equation}
t_{5}=t_{4}-\frac{1}{3} (t_{1} -t_{2} ), \\
\end{equation}
and the interior node  $t_{10}$ as,
\begin{equation}
t_{10} =\frac{1}{12} (t_{1} +t_{2} +4t_{3} +3t_{4} +3t_{5} )
\end{equation} 

The following algorithm depicts the overall discrete finite element computational procedure to obtain the numerical solution of the BVP. This algorithm outlines the steps for numerically solving the coupled system of a linear and a quasilinear partial differential equation using the finite element method, incorporating Picard's iteration for the nonlinear component.

\begin{algorithm}
\caption{Finite Element Method for Coupled Linear and Quasilinear PDEs}
\label{alg:coupled_fem}
\begin{algorithmic}[1]
\Require Model parameters ($\kappa$, $\xi$), Dirichlet boundary data for $\theta$ and $w$, computational domain $\Omega$, mesh refinement parameter $h$, convergence tolerance $\epsilon_{tol}$, maximum iterations $N_{max}$.
\Ensure Discrete solutions $\theta_h$ and $w_h$.

\State Generate a finite element mesh $\mathcal{T}_h$ for the domain $\Omega$.
\State Construct the discrete approximation spaces $\widehat{V}_h$ and $\widehat{V}_{h,0}$ based on the chosen shape functions (e.g., cubic-order Lagrange interpolants).
\State Apply the given Dirichlet boundary conditions for both $\theta$ and $w$ to constrain the solution spaces.

\Statex \Comment{Step 1: Solve the linear PDE for $\theta$}
\State Assemble the global stiffness matrix $A_{\theta}$ and load vector $F_{\theta}$ corresponding to the discrete weak formulation for $\theta$:
\Statex \quad $\sum_{K \in \mathcal{T}_h} \int_K \kappa \, \nabla \theta_h \cdot \nabla \varphi_h \,d\bfa{x} = 0, \quad \forall \varphi_h \in \widehat{V}_{h, \, 0}$.
\State Solve the resulting linear system $A_{\theta} \Theta = F_{\theta}$ to obtain the discrete solution $\theta_h \in \widehat{V}_h$.

\Statex \Comment{Step 2: Solve the quasilinear PDE for $w$ using Picard's Iteration}
\State Initialize an appropriate initial guess $w^0_h \in \widehat{V}_h$. This can be obtained by solving a simplified linear version of the $w$-equation (e.g., with $\Phi(w^n_h)$ set to a constant) or a zero field.
\State Set iteration counter $n=0$.
\Repeat
    \State Use the current iterate $w^n_h$ to evaluate the spatially varying coefficient $\Phi(w^n_h)$.
    \State Assemble the global stiffness matrix $A_{w}^{(n)}$ and load vector $F_{w}^{(n)}$ for the $(n+1)^{th}$ iteration of $w$:
    \Statex \quad $\sum_{K \in \mathcal{T}_h} \int_K \Phi(w^n_h) \; \nabla w^{n+1}_h \cdot \nabla \psi_h \; d\bfa{x} = - \sum_{K \in \mathcal{T}_h} \int_K \xi \, \theta_h \, \psi_h \; d\bfa{x}, \quad \forall \psi_h \in \widehat{V}_{h, \, 0}$.
    \State Solve the resulting linear system $A_{w}^{(n)} W^{n+1} = F_{w}^{(n)}$ to obtain the new iterate $w^{n+1}_h \in \widehat{V}_h$.
    \State Calculate the relative error: $Error = \frac{\|w^{n+1}_h - w^n_h\|}{\|w^{n+1}_h\|}$.
    \State Update $w^n_h \leftarrow w^{n+1}_h$.
    \State Increment $n \leftarrow n+1$.
\Until {$Error < \epsilon_{tol}$ or $n \geq N_{max}$}

\State The converged solutions are $\theta_h$ and $w_h = w^n_h$.
\end{algorithmic}
\end{algorithm}

\pagebreak 

\section{Numerical Examples and Validation}\label{num_exp}
To rigorously assess the computational efficiency and accuracy of the proposed numerical method, we present a series of three distinct numerical examples. The first example involves a nonlinear boundary value problem (BVP) for which an exact manufactured solution is available. This particular setup is crucial for performing a robust convergence analysis. The second example addresses a nonlinear BVP specifically defined on a V-notch domain, which introduces geometric complexities. Finally, the third and most challenging example extends the V-notch problem to include multiple inclusions, further complicating the domain geometry. A significant challenge in this third scenario is the generation of an appropriate finite element mesh in the vicinity of these intricate curved inclusions. Our proposed method effectively mitigates this difficulty by requiring a remarkably small number of elements to discretize the regions surrounding the inclusions. Moreover, a key advantage of utilizing curved triangular elements in our approach is the complete elimination of discretization errors associated with the exact representation of these curved boundaries.

\subsection{Square Domain with a Manufactured Solution}
\begin{figure}[H]
\centering
\begin{tikzpicture}
\draw (0,0) -- (3,0) -- (3,3) -- (0,3) -- (0,0);
\node at (-0.3,-0.25) {$(0,0)$};
\node at (3.3,3.2) {$(1,1)$};
\node at (-0.4, 1.5) {$\Gamma_{D_{1}}$};
\node at (3.4, 1.5) {$\Gamma_{D_{2}}$};
\node at (1.5, -0.25) {$\Gamma_{D_{3}}$};
\node at (1.5, 3.2) {$\Gamma_{D_{4}}$};
\end{tikzpicture}
\caption{Square Domain}
\label{fig_square}
\end{figure}
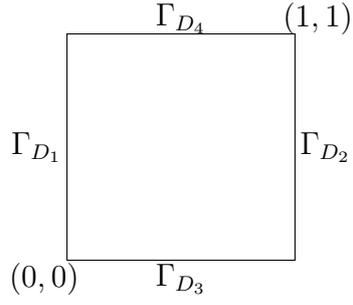

In this initial example, we analyze a nonlinear BVP defined over a square computational domain, denoted as $\Omega$. For validation purposes, we specifically choose the exact analytical solution for the field variable $w$ to be $w = \sin x \sin y$. This manufactured solution allows us to precisely determine the right-hand side term, $f$, by substituting $w$ directly into the strong form of the governing partial differential equation. Similarly, the Dirichlet boundary values are consistently prescribed using this known manufactured solution. This controlled environment is instrumental for conducting a comprehensive convergence analysis using the standard $L_2$ error measure, which quantifies the difference between our numerical approximation and the known exact solution. The precise layout of the computational domain is illustrated in Figure~\ref{fig_square}, and a detailed breakdown of the applied boundary conditions for this example is provided in Table~\ref{man_bc}.

\begin{equation}
\nabla . \left[ \frac{\nabla w}{1+ \nabla w} \right]  = - f(x,y) \quad in \quad \Omega
\end{equation}
with the boundary conditions in Table \ref{man_bc},
\begin{table}[H]
\centering
\caption{Boundary conditions for square domain}
\begin{tabular}{p{0.9 in} p{0.7 in}}\hline
 Boundary & Values  \\ \hline
 $\Gamma_{D_{1}}$ & 0\\
 $\Gamma_{D_{2}}$ & $0.8414709848 \sin y$\\
 $\Gamma_{D_{3}}$ & 0\\
 $\Gamma_{D_{4}}$ & $0.8414709848 \sin x$\\
\end{tabular}
 \label{man_bc}
\end{table}
The domain $\Omega$ is a square domain as shown in Figure \ref{fig_square}.

\begin{figure}[H]
 \centering
 \noindent
\includegraphics[scale=0.5]{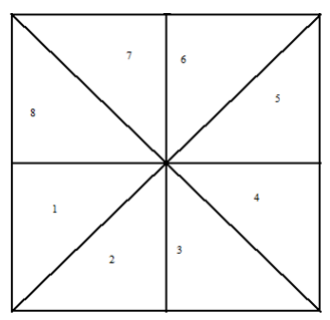}
 \caption{Manual Discretization of square domain into 8 elements}
 \label{fig_dis_square}
\end{figure}

Discretization of domain in Figure \ref{fig_square} into 8 elements manually as in Figure \ref{fig_dis_square}

\begin{table}[H]
\centering
\caption{Details of domain Figure \ref{fig_square} discretization}
\begin{tabular}{c c c}\hline
Elements & DOF & Boundary Points\\ \hline
8 & 49 & 24 \\ \hline
\end{tabular}
\label{tab_square}
\end{table}

\begin{table}[H]
\centering
\caption{Maximum $E_{a}$, $E_{r}$ and $L_2$-norm for 8 elements cubic-order triangular elements of Square domain}
\begin{tabular}{|c|c|c|c|} \hline
Iteartion Level & $E_{a}$ & $E_{r}$ & $L_{2}$ Norm\\ \hline
$0^{th}$ & 9.4213  $\times$ $10 ^{-5}$ & 0.34\% & 3.1280  $\times$ $10^{-4}$\\ \hline
$1^{st}$ & 1.1511  $\times$ $10 ^{-4}$ & 0.42\% & 3.2011  $\times$ $10^{-4}$\\ \hline
$2^{nd}$ & 1.1593  $\times$ $10 ^{-4}$ & 0.42\% & 3.2320  $\times$ $10^{-4}$\\ \hline
\end{tabular}
\label{manu_res}
\end{table}

\begin{figure}[H]
    \centering
    \begin{subfigure}[b]{0.48\textwidth} 
        \centering
        \includegraphics[width=\textwidth]{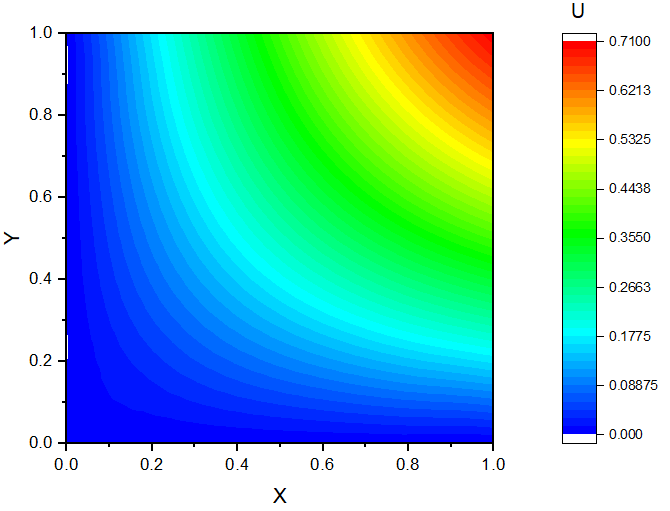}
        \caption{Contour Plot of Exact Solution}
        \label{fig:exact_solution}
    \end{subfigure}
    \hfill 
    \begin{subfigure}[b]{0.48\textwidth} 
        \centering
        \includegraphics[width=\textwidth]{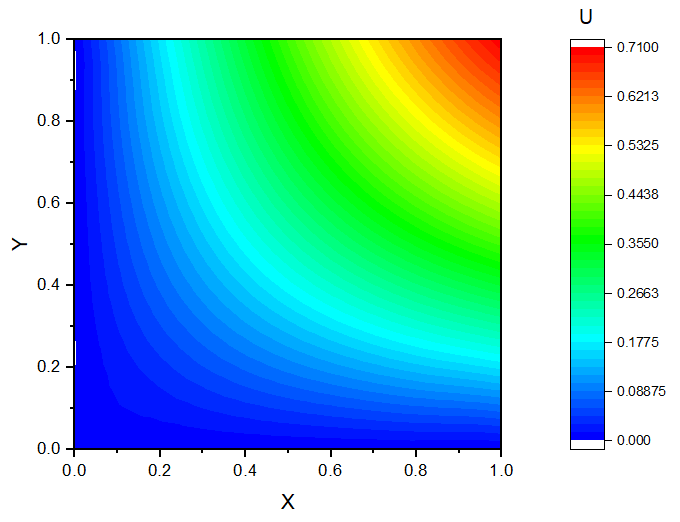}
        \caption{Contour Plot of FEM Solution}
        \label{fig:fem_solution}
    \end{subfigure}
    \caption{Comparison of Temperature Contour Plots: Exact vs. Finite Element Method (FEM) Solutions}
    \label{fig_contour_manufsoln}
\end{figure}

Meshing is a critical preliminary step in the Finite Element Method (FEM), particularly when dealing with complex, arbitrarily curved domains. The quality of the mesh directly impacts the accuracy and efficiency of the subsequent finite element analysis.

\subsection*{Proposed Meshing Scheme}
Our proposed meshing technique builds upon the robust DistMesh2D MATLAB code, originally developed by \cite{persson2004simple} for generating high-quality linear triangular elements. This method is favored for its simplicity and the excellent mesh characteristics it produces. The geometrical description within the DistMesh2D framework is remarkably concise, requiring only a few key inputs:
\begin{itemize}
    \item Signed distance function ($fd$): Defines the boundary of the domain.
    \item Edge length function ($fh$): Controls the desired element size across the domain.
    \item Initial mesh size ($h0$): Sets a global target for element size.
    \item Bounding box ($box$): Specifies the overall spatial extent of the domain.
    \item Fixed nodal positions ($fixed$): Allows for user-defined node placement at specific locations.
    \item Variable arguments ($varargin$): Provides flexibility for passing additional parameters to $fd$ and $fh$.
\end{itemize}

\subsection*{Extension to Curved Elements}
To accommodate curved triangular elements, our scheme incorporates the nodal relationship derived from subparametric transformations, a technique pioneered by \cite{rathod2008use}. This integration allows for the accurate representation of complex geometries using higher-order elements, which are essential for capturing intricate stress fields and behaviors in curved domains. Previous work by \cite{shylaja2019finite, shylaja2021improved, shylaja2021new, shylaja2021two, shylaja2025efficient} has successfully demonstrated the application of similar meshing approaches to various two-dimensional geometries, yielding high-quality meshes for higher-order curved triangular elements.

\subsection{Square Domain}
We're focusing on the numerical solution of the partial differential equations (PDEs) defined in Equation \eqref{f1a}. These PDEs are subject to specific boundary conditions detailed in Table \ref{ex1_bc}. The physical space, or computational domain, for this problem is a simple square domain, as visually represented in Figure \ref{fig_square}. To solve Equation \eqref{f1a} under the given boundary conditions, we discretize this square domain. For this particular analysis, the domain shown in Figure \ref{fig_square} is divided into eight discrete elements, as illustrated in Figure \ref{fig_dis_square}. This discretization process is a fundamental step in numerical methods like the Finite Element Method (FEM), transforming a continuous problem into a solvable discrete system.

Upon solving the discretized system, we can visualize the resulting solution. Figure \ref{fig_con_square_temp} presents a contour plot that clearly depicts the temperature profile across the square domain. This contour plot provides valuable insight into how the temperature varies spatially, fulfilling the conditions set by the PDEs and their associated boundary conditions.
\begin{table}[H]
\centering
\caption{Boundary Conditions for square domain}
\begin{tabular}{p{0.9 in} p{0.7 in} p{0.7 in}}\hline
 Boundary & Values \eqref{f1_1} & Values \eqref{f1_2}\\ \hline
 $\Gamma_{D_{1}}$ & 0 & 1\\
 $\Gamma_{D_{2}}$ & 0 & 0\\
 $\Gamma_{D_{3}}$ & $x(1-x)$ & $1-x$ \\
 $\Gamma_{D_{4}}$ & 0 & $1-x$\\
\end{tabular}
 \label{ex1_bc}
\end{table}

\begin{figure}[H]
    \centering
    \includegraphics[width=0.5\linewidth]{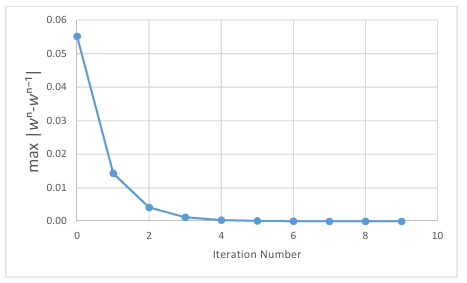}
    \caption{Plot showing the behavior of $max|w^n-w^{(n-1)}|$ at every Picard’s iteration for Square domain}
    \label{fig_picard_iteration_square}
\end{figure}

To solve Equation \eqref{f1_2}, we apply the boundary conditions detailed in Table \ref{ex1_bc}. The computational domain, as depicted in Figure \ref{fig_square}, is discretized into eight finite elements, a configuration visually represented in Figure \ref{fig_dis_square}. This discretization is a crucial step for numerical analysis, transforming the continuous problem into a system solvable by iterative methods. The convergence of our iterative solution process is illustrated in Figure \ref{fig_picard_iteration_square}. This plot shows the maximum absolute difference between two consecutive iterations, providing a clear indicator of how quickly the solution converges to a stable state. A decreasing trend in this plot signifies the successful iterative refinement of the solution.

Finally, after achieving convergence, we visualize a key mechanical output: the Airy stress profile. Figure \ref{fig_con_square} presents a contour plot of this profile, offering a comprehensive view of the stress distribution across the square domain. This visualization is essential for understanding the mechanical behavior of the material under the applied loads and boundary conditions.

\begin{figure}[H]
    \centering
    \begin{subfigure}[b]{0.48\textwidth} 
        \centering
        \includegraphics[width=\textwidth]{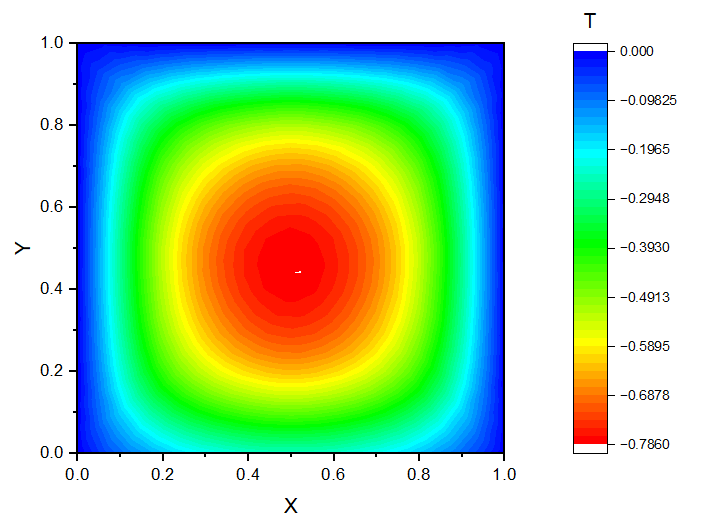}
        \caption{Temperature Profile}
        \label{fig:square_temp_profile}
    \end{subfigure}
    \hfill 
    \begin{subfigure}[b]{0.48\textwidth} 
        \centering
        \includegraphics[width=\textwidth]{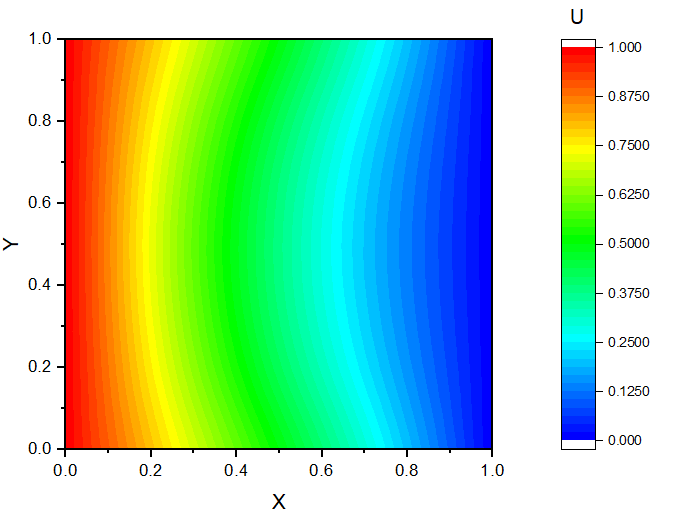}
        \caption{Airy Stress Profile}
        \label{fig:square_airy_stress_profile}
    \end{subfigure}
    \caption{Contour Plots for FEM Solution of a Square Domain: (a) Temperature Profile and (b) Airy Stress Profile}
    \label{fig_square_combined_contours}
\end{figure}

\subsection{V-notch Domain}
\begin{figure}[H]
 \centering
 \noindent\begin{subfigure}[b]{0.4\textwidth}
  \includegraphics[width=1\textwidth]{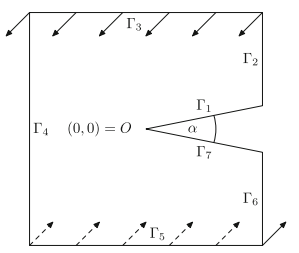}
 \end{subfigure}%
 \caption{V-notch Domain}
 \label{fig_vnotch}
\end{figure}

We're focusing on numerically solving the partial differential equations (PDEs) presented in Equation \eqref{f1a}. These PDEs govern the physical phenomenon under investigation and are subject to the specific boundary conditions meticulously outlined in Table \ref{ex2_bc}. The computational domain for this problem is particularly interesting: a V-notch geometry, as clearly illustrated in Figure \ref{fig_vnotch}. This type of domain often introduces stress concentrations or unique field behaviors, making its analysis crucial.

\begin{figure}[H]
\centering
\noindent\begin{subfigure}[b]{0.4\textwidth}
\includegraphics[width=1\textwidth]{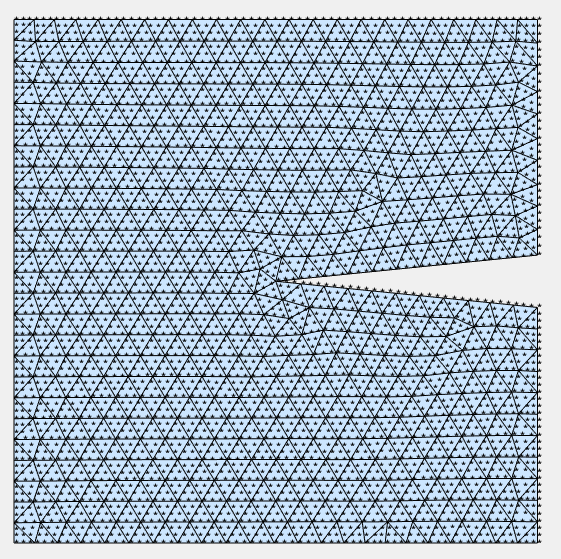}
\end{subfigure}%
\caption{MATLAB discretization of V-notch domain with edge length $h=0.045$}
\label{fig_dis_vnotch}
\end{figure}

To solve Equation \eqref{f1a} with the given boundary conditions, we employed a discretization approach using a MATLAB code. The resulting mesh, which divides the V-notch domain into discrete elements, is visually represented in Figure \ref{fig_dis_vnotch}. This meshing process is a foundational step in transforming the continuous PDE problem into a solvable system for numerical methods. Following the solution of the discretized system, the temperature profile across the V-notch domain is effectively visualized through a contour plot in Figure \ref{fig_con_vnotch_temp}. This plot provides valuable insight into how temperature is distributed throughout this complex geometry, allowing us to understand the thermal behavior under the specified conditions.

\begin{table}[H]
\centering
\caption{Boundary conditions for V-notch domain}
\begin{tabular}{p{0.9 in} p{0.7 in} p{0.7 in}}\hline
 Boundary & Values \ref{f1_1} & Values \ref{f1_2}\\ \hline
 $\Gamma_{1}$ & 0 & 0\\
 $\Gamma_{2}$ & 0 & 0\\
 $\Gamma_{3}$ & 0 & $1-x$\\
 $\Gamma_{4}$ & 0 & 1\\
 $\Gamma_{5}$ & $x(1-x)$ & $1-x$ \\
 $\Gamma_{6}$ & 0 & 0\\
 $\Gamma_{7}$ & 0 & 0\\
\end{tabular}
\label{ex2_bc}
\end{table}

\begin{table}[H]
\centering
\caption{Details of domain Fig.\ref{fig_dis_vnotch} discretization}
\begin{tabular}{p{0.6 in} p{0.5 in} p{1.25 in}}\hline
Elements & DOF & Boundary Points\\ \hline
1038 & 4843 & 342 \\ \hline
\end{tabular}
\label{tab_vnotch}
\end{table}

\begin{figure}[H]
    \centering
    \includegraphics[width=0.5\linewidth]{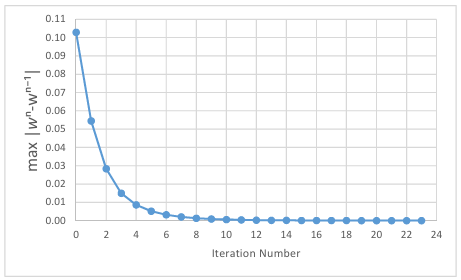}
    \caption{Plot showing the behavior of $max|w^n-w^{(n-1)}|$ at every Picard’s iteration for V-notch domain}
    \label{fig_picard_iteration_vnotch}
\end{figure}

To solve Equation \eqref{f1_2}, we applied the boundary conditions meticulously outlined in Table \ref{ex2_bc}. The computational domain, a V-notch geometry shown in Figure \ref{fig_vnotch}, was discretized for numerical analysis. This discretization involved dividing the domain into a total of 1038 elements, as visually represented by the mesh in Figure \ref{fig_dis_vnotch}. This high number of elements allows for a detailed and accurate representation of the complex stress fields within the notch. The convergence behavior of our iterative solution process is clearly illustrated in Figure \ref{fig_picard_iteration_vnotch}. This plot tracks the maximum absolute difference between two consecutive iterations, serving as a crucial indicator of the solution's stability and accuracy. A rapid decrease in this value signifies efficient convergence towards the final solution. Finally, after achieving convergence, we generated a contour plot to visualize the Airy stress profile across the V-notch domain. Figure \ref{fig_con_vnotch} provides this comprehensive visual, which is essential for understanding the stress distribution and identifying potential areas of stress concentration within this challenging geometry.

\begin{figure}[H]
    \centering
    \begin{subfigure}[b]{0.48\textwidth} 
        \centering
        \includegraphics[width=\textwidth]{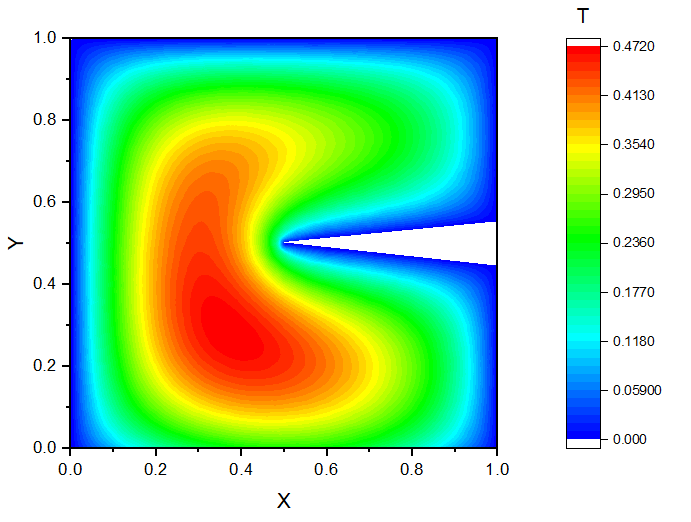}
        \caption{Contour Plot of temperature profile}
        \label{fig:picard_vnotch_sub}
    \end{subfigure}
    \hfill 
    \begin{subfigure}[b]{0.48\textwidth} 
        \centering
        \includegraphics[width=\textwidth]{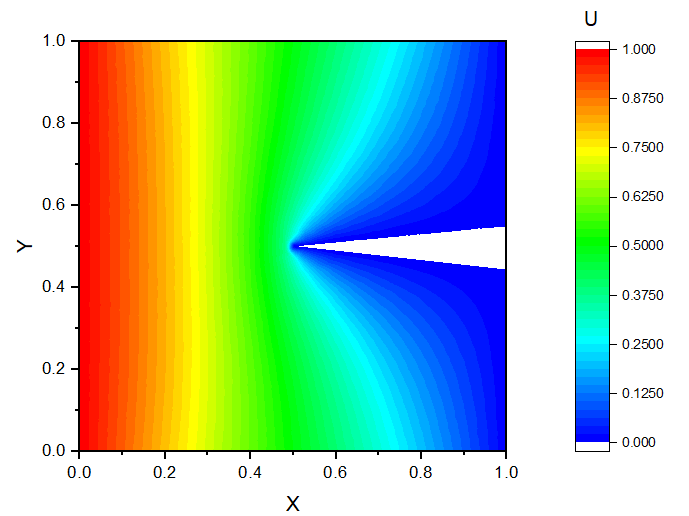}
        \caption{Contour Plot of Airy Stress Profile}
        \label{fig:airy_vnotch_sub}
    \end{subfigure}
    \caption{Numerical Results for V-Notch Domain: (a) Temperature profile and (b) Airy Stress Contour Plot}
    \label{fig_vnotch_combined_results}
\end{figure}

\subsection{V-notch with inclusions}
\begin{figure}[H]
 \centering
 \noindent\begin{subfigure}[b]{0.4\textwidth}
  \includegraphics[width=1\textwidth]{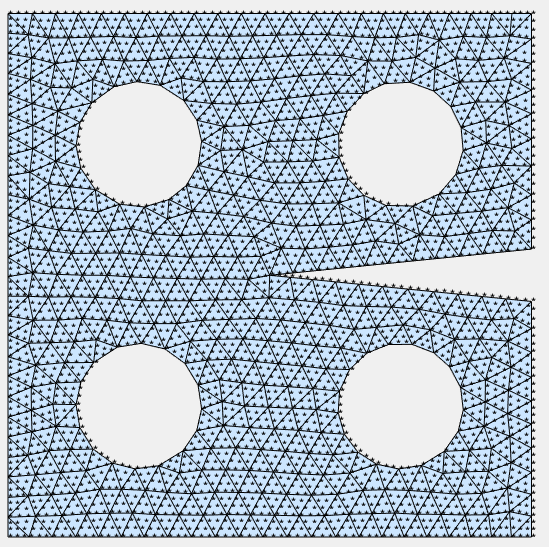}
 \end{subfigure}%
 \caption{MATLAB discretization of V-notch with inclusions domain with edge length $h=0.044$}
 \label{fig_vnotchwtinclusion}
\end{figure}

We turn our attention to the numerical solution of the partial differential equations (PDEs), as described by Equation \eqref{f1a}. For this particular scenario, the system is governed by the specific boundary conditions meticulously detailed in Table \ref{ex3_bc}. The computational domain for this analysis presents an increased level of geometrical complexity: a V-notch incorporating an inclusion. This configuration introduces heterogeneities that significantly influence the field behavior. The entire domain, including the V-notch and the embedded inclusion, was meticulously discretized using a custom MATLAB code. The resulting finite element mesh, which effectively represents the complex geometry, is clearly depicted in Figure \ref{fig_vnotchwtinclusion}. This intricate meshing is crucial for accurately capturing localized phenomena near the inclusion and the notch tip. We obtain the temperature profile across this complex domain after solving the discretized system. Figure \ref{fig_con_vnotchwtinclusions_temp} presents a contour plot that vividly illustrates this temperature distribution. This visualization is instrumental in understanding the thermal behavior and identifying any hot or cold spots influenced by the presence of both the V-notch and the inclusion.
\begin{table}[H]
\centering
\caption{Boundary conditions for V-notch domain}
\begin{tabular}{p{0.9 in} p{0.7 in} p{0.7 in}}\hline
 Boundary & Values \ref{f1_1} & Values \ref{f1_2}\\ \hline
 $\Gamma_{1}$ & 0 & 0\\
 $\Gamma_{2}$ & 0 & 0\\
 $\Gamma_{3}$ & 0 & $1-x$\\
 $\Gamma_{4}$ & 0 & 1\\
 $\Gamma_{5}$ & $x(1-x)$ & $1-x$ \\
 $\Gamma_{6}$ & 0 & 0\\
 $\Gamma_{7}$ & 0 & 0\\
\end{tabular}
 \label{ex3_bc}
\end{table}

\begin{table}[H]
\centering
\caption{Details of domain Figure \ref{fig_vnotchwtinclusion} discretization}
\begin{tabular}{p{0.6 in} p{0.5 in} p{1.25 in}}\hline
Elements & DOF & Boundary Points\\ \hline
836 & 4026 & 339 \\ \hline
\end{tabular}
\label{tab_vnotchwtincl}
\end{table}

\begin{figure}[H]
    \centering
    \includegraphics[width=0.5\linewidth]{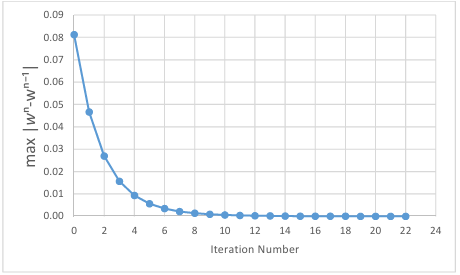}
    \caption{Plot showing the behavior of $max|w^n-w^{(n-1)}|$ at every Picard’s iteration for V-notch domain with inclusion}
    \label{fig_picard_iteration_vnotch_incl}
\end{figure}

\begin{figure}[H]
    \centering
    \begin{subfigure}[b]{0.48\textwidth} 
        \centering
        \includegraphics[width=\textwidth]{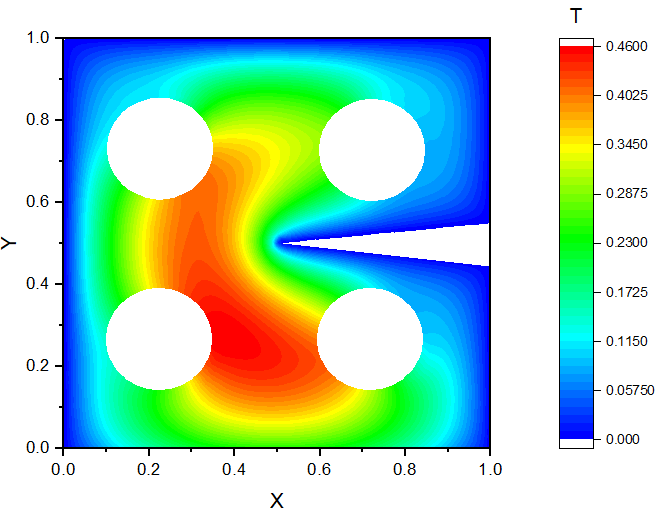}
        \caption{Temperature Profile}
        \label{fig:vnotch_incl_temp}
    \end{subfigure}
    \hfill 
    \begin{subfigure}[b]{0.48\textwidth} 
        \centering
        \includegraphics[width=\textwidth]{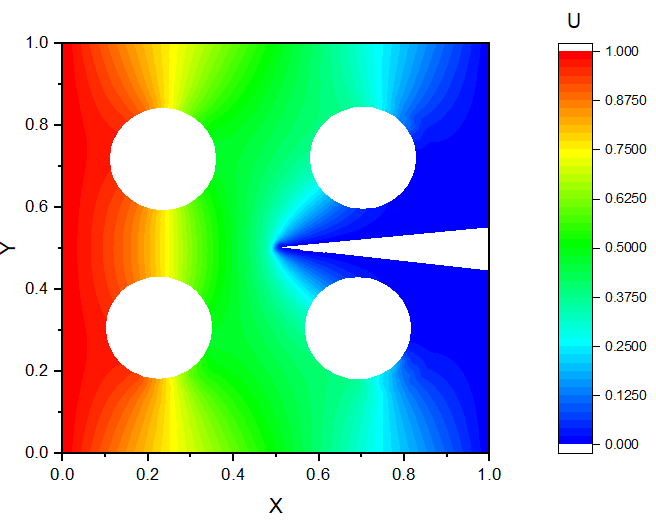}
        \caption{Airy Stress Profile}
        \label{fig:vnotch_incl_airy}
    \end{subfigure}
    \caption{Contour Plots for FEM Solution of a V-Notch Domain with Inclusion: (a) Temperature Profile and (b) Airy Stress Profile}
    \label{fig_vnotch_inclusion_combined_contours}
\end{figure}

\section{Conclusion}\label{conclusion}

This paper has presented a comprehensive finite element methodology for simulating the intricate behavior of a coupled thermoelastic system, specifically focusing on strain-limiting elastic materials. The proposed framework effectively addresses the challenges posed by an algebraically nonlinear and geometrically linear material response, particularly within domains featuring complex geometries such as V-notches and circular inclusions. Our formulation models this physical scenario as a coupled system of partial differential equations: a linear equation governing the temperature field and a quasilinear equation describing the mechanical deformation, derived from a specialized constitutive theory detailed in \cite{rajagopal2007elasticity,rajagopal2011modeling}. The inherent complexity of this coupled problem, particularly the quasilinear nature of the mechanical response, renders analytical closed-form solutions unattainable, thereby underscoring the critical need for robust numerical approaches.

A central innovation of this work lies in the strategic deployment of a novel class of cubic-order isoparametric curved triangular elements. These elements are uniquely capable of precisely representing the exact geometry of curved boundaries, such as those encountered around inclusions, thereby entirely circumventing any geometric discretization errors typically associated with standard straight-edged elements. Their seamless integration with conventional triangular meshes in the bulk material ensures a cohesive and accurate domain discretization. The numerical approximation of the primitive displacement variable benefits significantly from the high-order interpolation capabilities of these elements. To manage the intrinsic nonlinearity, a standard Picard iterative linearization scheme was successfully implemented, demonstrating rapid convergence to the solution across all test cases. The robustness and accuracy of our method were rigorously validated through its application to three distinct boundary value problems, including a critical comparison against a manufactured solution, where excellent agreement was consistently observed. This validation confirmed the enhanced accuracy of the approximation for the quasilinear partial differential equation.

In summary, the sophisticated methodology detailed herein offers an exemplary and vital tool with broad applicability. The developed approach for simulating elastic materials under combined thermal and mechanical loading is highly adaptable to a myriad of advanced computational mechanics problems. This includes, but is not limited to, the simulation of quasi-static crack propagation in both algebraically nonlinear elastic materials \cite{yoon2021quasi,lee2022finite} and materials whose parameters exhibit density dependence \cite{HCYSMMDDB2024}. Furthermore, the framework holds significant promise for investigations into multiscale phenomena \cite{MVSMM2023,vasilyeva2023} and complex multi-physics challenges \cite{SMMDDB2023}. A natural and important extension of this current work involves conducting a comprehensive convergence analysis of the entire coupled methodology, building upon findings similar to those presented in \cite{manohar2024hp}. A highly promising and compelling direction for future research involves the extension of the advanced numerical methodology presented in this paper, particularly its capabilities in handling quasilinear systems and complex geometries with curved elements. This extension would be invaluable for simulating quasi-static crack propagation in elastic materials, enabling a deeper understanding of fracture mechanics under various loading scenarios \cite{manohar2025adaptive, manohar2025convergence, fernando2025xi}. Moreover, adapting this robust framework to address full elasticity problems involving orthotropic materials represents another significant avenue for future investigation. This would considerably broaden the applicability of our method to a wider class of engineering materials that exhibit anisotropic mechanical properties, such as composites, where the material response varies with direction \cite{ghosh2025finite,ghosh2025computational}. Such developments would further solidify the versatility and practical utility of our finite element approach.

\bibliographystyle{plain}
\bibliography{coupling_reference}
\end{document}